\documentclass{article}
\usepackage[utf8x]{inputenc}

\usepackage{a4wide}
\usepackage{amscd}
\usepackage{amsmath}
\usepackage{amssymb}
\usepackage{amsthm}
\usepackage{array}
\usepackage{color}
\usepackage{epsfig}
\usepackage{graphics}
\usepackage{graphicx}
\usepackage[colorlinks=true,urlcolor=blue,linkcolor=blue,plainpages=false,pdfpagelabels
]{hyperref}
\usepackage{ifthen}
\usepackage{tabularx}
\usepackage{url}
\usepackage[all,cmtip]{xy}

\newcounter{uebno}

\newsavebox{\uebung}

\newcommand{\solutionsem}[1]{} 

\newtheorem{theorem}{Theorem}[section]

\newtheorem{definition}[theorem]{Definition}
\newtheorem{proposition}[theorem]{Proposition}
\newtheorem{corollary}[theorem]{Corollary}
\newtheorem{lemma}[theorem]{Lemma}
\newtheorem{fact}[theorem]{Remark}
\newtheorem{exemplu}[theorem]{Example}
\newtheorem{exercise}{Exercise}
\newtheorem{notation}[theorem]{Notation}

\newcommand{\bdfn}{\begin{definition}}
\newcommand{\edfn}{\end{definition}}
\newcommand{\bthm}{\begin{theorem}}
\newcommand{\ethm}{\end{theorem}}
\newcommand{\bprop}{\begin{proposition}}
\newcommand{\eprop}{\end{proposition}}
\newcommand{\bcor}{\begin{corollary}}
\newcommand{\ecor}{\end{corollary}}
\newcommand{\blem}{\begin{lemma}}
\newcommand{\elem}{\end{lemma}}
\newcommand{\bfact}{\begin{fact}}
\newcommand{\efact}{\end{fact}}
\newcommand{\bex}{\begin{exemplu}\begin{rm}}
\newcommand{\eex}{\end{rm}\end{exemplu}}
\newcommand{\bxc}{\begin{exercise}}
\newcommand{\exc}{\end{exercise}}
\newcommand{\bntn}{\begin{notation}}
\newcommand{\entn}{\end{notation}}

\newcommand{\be}{\begin{enumerate}}
\newcommand{\ee}{\end{enumerate}}
\newcommand{\bce}{\begin{center}}
\newcommand{\ece}{\end{center}}
\newcommand{\bi}{\begin{itemize}}
\newcommand{\ei}{\end{itemize}}
\newcommand{\bt}{\begin{tabular}}
\newcommand{\et}{\end{tabular}}
\newcommand{\beq}{\begin{equation}}
\newcommand{\eeq}{\end{equation}}
\newcommand{\ba}{\begin{array}} 
\newcommand{\ea}{\end{array}}
\newcommand {\bea} {\begin{eqnarray}}
\newcommand {\eea} {\end {eqnarray}}
\newcommand {\bua} {\begin{eqnarray*}}
\newcommand {\eua} {\end {eqnarray*}}

\newcommand{\se}{\subseteq}

\newcounter{ct}

\def\R{{\mathbb R}}
\def\N{{\mathbb N}}

\def\R{{\mathbb R}}

\newcommand{\eps}{\varepsilon}

\newcommand{{\Fdo}}{F}

\title{A note on the Mann iteration for $k$-strict pseudocontractions in Banach spaces}
\author{Andrei Sipo\c s${}^{a,b}$\\[0.2cm]
\footnotesize ${}^a$Faculty of Mathematics and Computer Science, University of Bucharest,\\
\footnotesize Academiei 14, 010014 Bucharest, Romania\\[0.1cm]
\footnotesize ${}^b$Simion Stoilow Institute of Mathematics of the Romanian Academy,\\
\footnotesize P. O. Box 1-764, 014700 Bucharest, Romania\\[0.1cm]
\footnotesize E-mail: Andrei.Sipos@imar.ro
}
\date{}

\begin{document}

\maketitle

\begin{abstract}
We show that a variant of a function defined by Cholamjiak and Suantai can be used to characterize 2-uniform smoothness. We then obtain greatly simplified proofs of a convergence theorem of Marino and Xu and of one of Zhou using a generalization of a lemma of Browder and Petryshyn and the aforementioned characterization. We also obtain more easily rates of asymptotic regularity corresponding to the studied iterations. Finally, we derive a way to relate two constants which are characteristic to 2-uniformly smooth spaces.


\end{abstract}

\section{Introduction}

The class of $k$-strict pseudocontractions was introduced by Browder and Petryshyn in \cite{BroPet67} for Hilbert spaces. If $H$ is a Hilbert space, $C \se H$ is a convex subset and $k \in [0,1)$, then a mapping $T: C \to H$ is called a {\bf $k$-strict pseudocontraction} if for all $x,y \in C$ we have that:
\beq
\|Tx-Ty\|^2 \leq \|x-y\|^2 + k\|(x-Tx) - (y-Ty)\|^2.\label{k-ps}
\eeq
If we set $k:=0$ in the above, we obtain the condition $\|Tx-Ty\| \leq \|x-y\|$, which states that the mapping $T$ is nonexpansive. The search of algorithms for finding fixed points of nonexpansive self-mappings of subsets of metric spaces belonging to various established classes has been a longstanding research program. In the sequel, we will restrict ourselves to classes of Banach spaces. We denote the unit sphere of a Banach space $E$ by $S(E)$. A Banach space $E$ is called {\bf smooth} if for any $u \in S(E)$, its norm is Gâteaux differentiable at $u$, i.e. for any $v \in S(E)$, the limit
$$\lim_{h \to 0} \frac{\|u+hv\|-\|u\|}{h}$$
exists. This condition has been proven to be equivalent to the fact that the normalized duality mapping of the space, $J: E \to 2^{E^*}$, is single-valued -- and we shall denote its unique section by $j: E \to E^*$. Therefore, for all $x \in E$, $j(x)(x)=\|x\|^2$ and $\|j(x)\|=\|x\|$. Moreover, $E$ has a {\bf Fréchet differentiable norm} if, in addition, the limit above is attained uniformly in the variable $v \in S(E)$ and it is {\bf uniformly smooth} (or has a {\bf uniformly Fréchet differentiable norm}) if the limit is attained uniformly in the pair of variables $(u,v) \in S(E) \times S(E)$. One can also define the {\bf modulus of smoothness} of $E$ to be the map $\rho_E: (0, \infty) \to \R$, defined, for all $\tau \in (0, \infty)$, by
$$\rho_E(\tau):= \sup\left\{\frac{\|u+\tau v\|+\|u-\tau v\|}{2} -1\ \middle|\ u,v\in S(E) \right\}.$$
It is known that a space $E$ is uniformly smooth iff 
$$ \lim_{\tau\to 0}\frac{\rho_E(\tau)}{\tau}=0.$$
This can happen, for example, if there are $c>0$ and $q>1$ such that for all $\tau$, $\rho_E(\tau) \leq c \tau^q$. In that case, $E$ is said to be {\bf $q$-uniformly smooth}.

We define now the {\bf modulus of convexity} of a space $E$ to be the map $\delta_E: [0,2] \to \R_+$, defined, for all $\eps \in [0,2]$, by
$$\delta_E(\eps):= \inf\left\{1-\frac{\|x+y\|}{2}\ \middle|\ x,y\in S(E), \|x-y\|\geq\eps \right\}.$$
A Banach space $E$ is called {\bf uniformly convex} iff for all $\eps>0$, $\delta_E(\eps)>0$ -- or, equivalently, if for any $\eps >0$ there is a $\delta>0$ such that for any $u,v \in S(E)$ such that $\|u-v\|\geq\eps$ we have that $\frac{\|u+v\|}{2}\leq 1-\delta$. If $\eta : \R \to \R$ is a function such that for any $\eps>0$, $\eta(\eps)$ is such a $\delta$, we call $\eta$ a {\bf valid modulus of uniform convexity} for the space. We remark that the modulus of convexity from above is, in a sense, the ``optimal'' valid modulus of uniform convexity.

The algorithms mentioned above are usually iterative in nature -- a typical example is the {\bf Mann iteration} associated to such a self-mapping $T: C \to C$, an initial point $x \in C$ and a sequence $(t_n)_{n \in \N} \subseteq (0,1)$, which is the sequence $(x_n)_{n \in \N} \subseteq C$ defined\footnote{We note that for the definition to make sense, $C$ should be presupposed to be convex.} by:
$$x_0:=x$$
$$x_{n+1}:=t_n Tx_n + (1-t_n) x_n$$

A result on the Mann iteration for nonexpansive mappings is the following theorem of Reich, which will be needed in the sequel.

\bthm[Reich (1979), {\cite[Theorem~2]{Rei79}}]\label{reich}
Let $E$ be a uniformly convex Banach space with a Fréchet differentiable norm, $C \se E$ a convex, closed set and $T:C \to C$. Suppose that $T$ is nonexpansive with $Fix(T)\neq\emptyset$. Let $x \in C$ and $(t_n)_{n \in \N} \subseteq (0,1)$ such that
$$\sum_{n=0}^\infty t_n(1-t_n) = \infty$$
Then the Mann iteration corresponding to $T$, $x$ and $(t_n)_{n \in \N}$ weakly converges to a fixed point of $T$.
\ethm

The class of $k$-strict pseudocontractions has been less readily amenable to classical iterative schemas. Still, in 2007, Marino and Xu proved\footnote{Marino and Xu use the notation $\alpha_n:=1-t_n \in (k,1)$, so that the condition becomes $\sum_{n=0}^\infty (\alpha_n-k)(1-\alpha_n) = \infty$. We have chosen to present it as above in order to maintain uniformity with Reich's approach.} the following theorem.

\bthm[Marino and Xu (2007), {\cite[Theorem~3.1]{MarXu07}}]\label{marxu}
Let $H$ be a Hilbert space, $C\se H$ a convex, closed set and $T: C \to C$. Let $k$ be in $(0,1)$ and suppose that $T$ is a $k$-strict pseudocontraction with $Fix(T)\neq\emptyset$. Let $x \in C$ and $(t_n)_{n \in \N} \subseteq (0,1-k)$ such that
$$\sum_{n=0}^\infty t_n(1-k-t_n) = \infty$$
Then the Mann iteration corresponding to $T$, $x$ and $(t_n)_{n \in \N}$ weakly converges to a fixed point of $T$.
\ethm

Marino and Xu asked in their paper whether this result can be generalized to uniformly convex Banach space with a Fréchet differentiable norm, in the same vein as Reich's. Since then, various authors have tried to solve this problem to some degree.

The first issue that arises is that of the proper generalization of $k$-strict pseudocontractions to the case of Banach spaces. The solution (used, for example, in \cite{ChoSua13,Zho14}) comes from the initial observation (\cite[Theorem 1.(2)]{BroPet67}) that condition \eqref{k-ps} is equivalent to the following:
$$\langle (x-Tx)-(y-Ty), x-y \rangle \geq \frac{1-k}{2} \|(x-Tx) - (y-Ty)\|^2.$$
Now, if $E$ is a smooth Banach space, as we said above, it admits a single-valued normalized duality mapping $j: E \to E^*$. Hence the natural extension of the condition above in this framework is the following one, which we will take as our official definition of $k$-strict pseudocontractions in Banach spaces:
$$j(x-y)((x-Tx)-(y-Ty)) \geq \frac{1-k}{2} \|(x-Tx) - (y-Ty)\|^2.$$

A significant advance in this direction was made by Zhou in 2014, who, using a variant\footnote{The original definition in \cite[Lemma 3.2]{ChoSua13} was: 
$$\beta^*_E(x,t):=\sup\left\{\middle|\frac{\|x+tv\|^2-\|x\|^2}t-2j(x)(v)\middle|\ \middle|\ v \in S(E)\right\}.$$
which would make Zhou's condition \eqref{betastar} unnecessarily stronger.} of a function previously defined by Cholamjiak and Suantai -- namely, for any Banach space $E$ with a Fréchet differentiable norm, one can define the function $\beta^*_E: E \times (0,\infty) \to \R$ by:
$$\beta^*_E(x,t):=\sup\left\{\frac{\|x+tv\|^2-\|x\|^2}t-2j(x)(v)\ \middle|\ v \in S(E)\right\},$$
for any $x \in E$ and $t \in (0, \infty)$ -- proved the following result, which is implicitly contained in \cite[Theorem 3.1]{Zho14}.

\bthm[Zhou (2014), \cite{Zho14}]\label{the-zhou}
Let $E$ be a uniformly convex Banach space which is also uniformly smooth, $C\se E$ a convex, closed set and $T: C \to C$. Let $d \in [1, \infty)$ be such that for any $x \in E$ and $t \in (0, \infty)$,
\beq
\beta^*_E(x,t)\leq dt.\label{betastar}
\eeq
Let $k$ be in $(0,1)$ and suppose that $T$ is a $k$-strict pseudocontraction with $Fix(T)\neq\emptyset$. Let $x \in C$ and $(t_n)_{n \in \N} \subseteq (0,\frac{1-k}{2d})$ such that
$$\sum_{n=0}^\infty t_n = \infty.$$
Then the Mann iteration corresponding to $T$, $x$ and $(t_n)_{n \in \N}$ weakly converges to a fixed point of $T$.
\ethm

Our aim is to show that condition \eqref{betastar} above is actually equivalent to $2$-uniform smoothness and then to give simpler and immediate proofs of both Theorem~\ref{marxu} and Theorem~\ref{the-zhou}.

\section{The main results}

By convention, we will set for any space $E$ and any $x \in E$, $\beta^*_E(x,0):=0$. We first note that for any Hilbert space $H$, any $x\in H$ and any $t \geq 0$, $\beta^*_H(x,t)=t$. Indeed, in that case, for any $x \in H$, $t>0$ and $v \in S(H)$, we have that:
$$\frac{\|x+tv\|^2-\|x\|^2}t-2j(x)(v) = \frac{\|x+tv\|^2-\|x\|^2 -2\langle x,tv \rangle}t = \frac{\|x+tv-x\|^2}t = t.$$

The characterization lemma is the following.

\blem\label{l1}
Let $E$ be a smooth Banach space. The following statements are equivalent:
\be
\item $E$ is $2$-uniformly smooth, i.e. there is a $c>0$ such that for all $\tau$, $\rho_E(\tau) \leq c \tau^2$;
\item there is a $d>0$ such that for all $x,y \in E$ we have that $\|x+y\|^2 \leq \|x\|^2 +2j(x)(y)+ d\|y\|^2$;
\item there is a $d>0$ such that for all $x \in E$ and $t \geq 0$ we have that $\beta^*_E(x,t)\leq dt$.
\ee
Moreover, the constants in (ii) and (iii) may be taken to be the same (and hence we have used the same designator).
\elem

\begin{proof}
The equivalence between (i) and (ii) is given by \cite[Corollary 1']{Xu91}.

Suppose now that (ii) holds. Let $x \in E$ and $t\geq 0$. If $t=0$, there is nothing to prove, so suppose $t>0$. Let $v$ be in $S(E)$. Then, by (ii), setting $y:=tv$, we get that:
$$\|x+tv\|^2 \leq \|x\|^2 +2tj(x)(v)+ dt^2,$$
so
$$\frac{\|x+tv\|^2-\|x\|^2}t-2j(x)(v) \leq dt,$$
from which, by taking the supremum, the conclusion follows.

Suppose now that (iii) holds. Let $x,y \in E$. Again, if $y=0$, there is nothing to prove. If $y \neq 0$, set $v:=\frac{1}{\|y\|}\cdot y$ and $t:=\|y\|$, so $tv=y$. Then we get that:
$$\frac{\|x+y\|^2-\|x\|^2}{\|y\|}-\frac{2}{\|y\|}\cdot j(x)(y) \leq \beta^*_E(x,\|y\|) \leq d\|y\|,$$
and by multiplying by $\|y\|$ we obtain our desired result.
\end{proof}

We have therefore established the equivalence of Zhou's condition with $2$-uniform smoothness, which was mentioned as a special case in \cite[p. 762]{Zho14}. We note that it is immediate that if $d$ satisfies the two equivalent conditions, and $d\leq d'$, then $d'$ also satisfies the condition. Hence we can always take $d \geq 1$. In addition to that, we will also sketch in Section~\ref{sec:dc} a way to compute the constant $d$ -- and that specific choice of $d$ for a given $c$ will be denoted by $d_c$ from now on (and its value is already greater than $1$).

The other ingredient of our result is the following generalization of \cite[Theorem 2]{BroPet67}. For a given self-mapping of a convex set, $T: C \to C$, and a $t \in (0,1)$, set $T_t:= tT+(1-t)id_C$ -- that is, for all $x \in C$, $T_tx = tTx + (1-t)x$. It is immediate that for all $t_1, t_2 \in (0,1)$, $(T_{t_1})_{t_2}=T_{t_1 \cdot t_2}$. Also note that, for any $t$, $T$ and $T_t$ have the same fixed points.

\blem\label{l2}
Let $E$ be a Banach space, $C \se E$ a convex subset and $d\geq 1$ such that for any $x \in E$ and $t \leq 0$, $\beta^*_E(x,t) \leq dt$. Let $k \in (0,1)$ and $T : C \to C$ a $k$-strict pseudocontraction.

Let $t \in (0,\frac{1-k}{d}]$. Then $T_t$ is nonexpansive. (In particular, $T_{\frac{1-k}{d}}$ is nonexpansive.)
\elem

\begin{proof}
Since $t \leq \frac{1-k}{d}$, we have that $dt-(1-k) \leq 0$.

Let $x,y \in E$. We have that:
\begin{align*}
 \|T_tx - T_ty\|^2 &= \|tTx + (1-t)x -tTy - (1-t)y\|^2 \\
 &= \|(x-y) + (-t)((x-Tx)-(y-Ty))\|^2 \\
 &\leq \|x-y\|^2 -2tj(x-y)((x-Tx)-(y-Ty)) + dt^2\|(x-Tx)-(y-Ty)\|^2\\
 &\leq \|x-y\|^2 -t(1-k)\|(x-Tx)-(y-Ty)\|^2  + dt^2\|(x-Tx)-(y-Ty)\|^2\\
 &= \|x-y\|^2 +t(dt - (1-k))\|(x-Tx)-(y-Ty)\|^2 \\
 &\leq \|x-y\|^2.\\
\end{align*}
\end{proof}

Now we can prove the main convergence result.

\bthm\label{main}
Let $E$ be a uniformly convex Banach space which is also $2$-uniformly smooth, $C\se E$ a convex, closed set and $T: C \to C$. Let, therefore, $d\geq 1$ be a constant satisfying conditions (ii) and (iii) from Lemma~\ref{l1} (if, for example, $\rho_E(\tau) \leq c \tau^2$, for all $\tau$, take $d:=d_c$). Let $k$ be in $(0,1)$ and suppose that $T$ is a $k$-strict pseudocontraction with $Fix(T)\neq\emptyset$. Let $x \in C$ and $(t_n)_{n \in \N} \subseteq (0,\frac{1-k}{d})$ such that
$$\sum_{n=0}^\infty t_n\left(\frac{1-k}{d}-t_n\right) = \infty$$
Then the Mann iteration corresponding to $T$, $x$ and $(t_n)_{n \in \N}$ weakly converges to a fixed point of $T$.
\ethm

\begin{proof}
By Lemma~\ref{l2}, we have that $T_{\frac{1-k}{d}}$ is nonexpansive. For every $n \geq 0$, set $t'_n:=t_n \cdot \frac{d}{1-k}$. Denote by $(x_n)_{n \in \N}$ the Mann iteration corresponding to $T$, $x$ and $(t_n)_{n \in \N}$. Let $n \geq 0$. We have that:

\begin{align*}
x_{n+1} &= t_nTx_n + (1-t_n)x_n \\
&= T_{t_n}x_n \\
&= T_{t'_n \cdot \frac{1-k}{d}} x_n \\
&= T_{t'_n} (T_{\frac{1-k}{d}} x_n) \\
&= t'_n T_{\frac{1-k}{d}} x_n + (1-t'_n) x_n .
\end{align*}
We have then, that $(x_n)_{n \in \N}$ is the Mann iteration corresponding to $T_{\frac{1-k}{d}}$, $x$ and $(t'_n)_{n \in \N}$. We seek to apply Theorem~\ref{reich}. For that we do the following verification:
$$\sum_{n=0}^\infty t'_n (1-t'_n) = \sum_{n=0}^\infty t_n \cdot \frac{d}{1-k} \left(1- t_n \cdot \frac{d}{1-k} \right) = \left(\frac{d}{1-k}\right)^2  \sum_{n=0}^\infty t_n\left(\frac{1-k}{d}-t_n\right) = \infty.$$
We therefore get that $(x_n)_{n \in \N}$ weakly converges to a fixed point of $T_{\frac{1-k}{d}}$, which is also a fixed point of $T$.
\end{proof}

\begin{proof}[Proof of Theorem \ref{the-zhou}]
Note that the hypothesis states that $(t_n)_{n \in \N} \subseteq (0,\frac{1-k}{2d})$. Then, for all $n$, $\frac{1-k}{d} - t_n \geq \frac{1-k}{2d}$, so:
$$\sum_{n=0}^\infty t_n\left(\frac{1-k}{d}-t_n\right) \geq \frac{1-k}{2d}\sum_{n=0}^\infty t_n = \infty.$$
We are therefore in the hypothesis of Theorem~\ref{main}.
\end{proof}

\begin{proof}[Proof of Theorem \ref{marxu}]
We have shown in the beginning of this section that one can take for a Hilbert space $d:=1$. The conclusion immediately follows.
\end{proof}

Note that for directly establishing the truth of Theorem \ref{marxu}, we would not have needed the full machinery of the lemmas, but only the result of \cite[Theorem 2]{BroPet67}. Also take into consideration, as a reviewer pointed out, that our results above can be considered special cases of the results of \cite{ColMar11}, specifically Lemma 3 and Theorem 12 (modulo different notations). Still, we feel that the above exposition, not being burdened by the extraneous details of the cyclic algorithm considered in \cite{ColMar11}, can shed light on the way in which earlier results on strict pseudocontractions, like the one in \cite{MarXu07}, can be deduced from older results like Reich's without rebuilding from the ground up the whole machinery on demiclosedness and monotonicity.

\section{Quantitative results}

Proof mining is a research program introduced by U. Kohlenbach in the 1990s (\cite{Koh08} is a comprehensive reference), which aims to obtain explicit quantitative information (witnesses and bounds) from proofs of an apparently ineffective nature. This paradigm in applied logic has successfully led so far to obtaining some previously unknown effective bounds, primarily in nonlinear analysis and ergodic theory. A large number of these are guaranteed to exist by a series of logical metatheorems which cover general classes of bounded or unbounded metric structures.

An example of such a piece of quantitative information is the following.

For any metric space $(X,d)$, any $T : X \to X$ and any $(x_n)_{n \in \N}$, we say that the sequence $(x_n)_{n \in \N}$ is {\bf $T$-asymptotically regular} if
$$\lim_{n \to \infty} d(x_n,Tx_n) = 0.$$
We call a {\bf rate of $T$-asymptotic regularity} for $(x_n)_{n \in \N}$ a rate of convergence for the sequence above, i.e. a function $h : (0, \infty) \to \N$ such that for any $\eps >0$ and any $n \geq h(\eps)$,
$$ d(x_n,Tx_n) \leq \eps.$$

Kohlenbach has computed, using proof mining techniques, in \cite[Theorem 3.4]{Koh03} that the Mann iteration of the conclusion of Theorem~\ref{reich} has a rate of $T$-asymptotic regularity of:
$$h^{(1)}_{b,\theta,\eta}(\eps) := \theta \left( \left\lceil \frac{3(b+1)}{2\eps \cdot \eta(\frac{\eps}{b+1})} \right\rceil \right),$$
where $\eta$ is a valid modulus of uniform convexity for the space, $b$ is a bound on the minimum distance between the initial point $x$ and a fixed point of $T$ and $\theta: \N \to \N$ is such that for any $N$,
$$\sum_{n=0}^{\theta(N)} t_n (1-t_n) \geq N.$$
Such a $\theta$ is called a {\bf rate of divergence} for the series.

Moreover, using the remarks from the statement of the theorem in \cite{Koh03} it can be shown that in the case of Hilbert spaces, which have the well-known modulus of uniform convexity of $\eta(\eps):=\frac{\eps^2}{8}$, one can simplify the above rate to:
\beq
h^{(2)}_{b,\theta}(\eps) := \theta \left( \left\lceil \frac{4(b+1)}{\eps^2 } \right\rceil \right). \label{hilb}
\eeq

The proof of Theorem~\ref{main} shows that if $\theta: \N \to \N$ is a rate of divergence for the series
$$\sum_{n=0}^\infty t_n\left(\frac{1-k}{d}-t_n\right),$$
then, since $\left \lceil \left( \frac{1-k}{d} \right) \right \rceil = 1$ (remember that we always took $d \geq 1$), the same $\theta: \N \to \N$ is also a rate of divergence for the series
$$\sum_{n=0}^\infty t'_n (1-t'_n).$$
In addition, an easy computation gives us that:
$$\|x_n - T_{\frac{1-k}{d}}x_n\|=\frac{1-k}{d}\|x_n-Tx_n\|.$$

We can, therefore, give the following result:
\bthm
A rate of $T$-asymptotic regularity for the sequence in Theorem~\ref{main} is given by:
$$h^{(3)}_{b,\theta,\eta}(\eps):=\theta\left(\left\lceil \frac{3(b+1)d}{2\eps(1-k) \cdot \eta\left(\frac{\eps(1-k)}{(b+1)d}\right)} \right\rceil \right),$$
where $b$, $\eta$ are as before and $\theta$ is a rate of divergence for the series
$$\sum_{n=0}^\infty t_n\left(\frac{1-k}{d}-t_n\right).$$
\ethm

\begin{proof}
Applying Kohlenbach's result and the remark on $\theta'$ from above, we get that for any $n \geq h^{(3)}_{b,\psi,\eta}(\eps)$,
$$\|x_n - T_{\frac{1-k}{d}}x_n\| \leq \frac{\eps(1-k)}{d}.$$
Now, by the computation before, we get that:
$$\|x_n-Tx_n\| \leq \frac{d}{1-k} \cdot \frac{\eps(1-k)}{d} = \eps.$$
\end{proof}

Applying the same treatment to the Hilbert rate from \eqref{hilb}, we get that a rate of $T$-asymptotic regularity in the case of Theorem~\ref{marxu} is:
$$h^{(4)}_{b,\theta}(\eps) := \theta \left(  \left\lceil \frac{4(b+1)}{(1-k)^2\eps^2 } \right\rceil \right),$$
a quadratic rate not unlike the $\theta \left( \left\lceil \frac{b^2}{\eps^2 } \right\rceil \right)$ one, obtained also with proof mining techniques in \cite{IvaLeu15}, for the same iteration.

\section{Computing the constant}\label{sec:dc}

We proceed to derive a formula for $d_c$. Let $E$ be a $2$-uniformly smooth Banach space and let, therefore, $c>0$ be such that $\rho_E(\tau) \leq c \tau^2$, for all $\tau$. We shall use Lindenstrauss's classical formula from \cite{Lin63}:
$$\rho_E(t) = \sup_{\eps\in [0,2]}\ \left(\frac{1}{2}\eps t - \delta_{E^*}(\eps) \right).$$
So for any suitable $\eps$ and $t$, we get that:
$$\delta_{E^*}(\eps) \geq \frac{1}{2}\eps t -ct^2.$$
The term on the right hand side is a quadratic function of $t$, which has as its maximum the value $\frac{\eps^2}{16c}$. So for all $\eps$, we have that:
$$\delta_{E^*}(\eps) \geq \frac{1}{16c} \cdot \eps^2.$$

Lemma 1 from \cite{Fig72} states that in this case, for any $x,y \in E^*$ with
$$\|x\|^2 + \|y\|^2 = 2$$
the following holds:
$$\frac{1}{2}\|x+y\| \leq 1- k_1 \left( \frac{1}{2} \|x-y\| \right)^2,$$
where $k_1$ is given by the minimum of $\frac{1}{16c}$ and $\alpha$, the constant $\alpha$ being such that
$$\sup_{t \in (0,1]} \frac{\sqrt{2-(1-t)^2}-1-t}{t^2}=-\alpha <0.$$
The function to be maximized here can be seen to be equal to $\frac{-2}{t+1+\sqrt{2-(1-t)^2}}$, which increases along with its obviously increasing denominator. So the original function attains its maximum at $t:=1$, the maximum being $\sqrt{2}-2$. The value of $\alpha$ is therefore $2-\sqrt{2}$.

We denote by $L_2(E^*)$ the space of all functions $f:[0,1]\to E^*$ such that:
$$\int_0^1 \|f(t)\|^2 dt < \infty.$$
By Proposition 1 and ``Added in proof'' of \cite{Fig72}, we get that for all $\eps \in [0,2]$,
$$\delta_{L_2(E^*)}(\eps) \geq k_2 \eps^2,$$
where $k_2$ is computed by
$$k_2:= 2^{-2} \cdot 2^{-1} \cdot \min(k_1,1) = \frac{\min\left(\frac{1}{16c},2-\sqrt{2}\right)}{8}.$$

By the statement and proof of Lemma 2.1 from \cite{PruSma87}, we get that for all $x,y \in E^*$ and all $t \in (0,1)$,
$$\|tx+(1-t)y\|^2 \leq t\|x\|^2 + (1-t)\|y\|^2 - k_2t(1-t)\|x-y\|^2,$$
which corresponds to equation (3.1) from \cite{Xu91}.

Then, using the proof of the implication (i) $\Rightarrow$ (ii) from \cite[Corollary 1]{Xu91}, we get that for all $x,y \in E^*$ and all $f \in J(x)$,
$$\|x+y\|^2 \geq \|x\|^2 + 2f(y) + k_2\|y\|^2.$$

Finally, from the first lines of the proof of \cite[Theorem 1']{Xu91}, we obtain that for all $x,y \in E$,
$$\|x+y\|^2 \leq \|x\|^2 + 2j(x)(y) + k_2^{-1}\|y\|^2.$$

We can therefore set
$$d_c := k_2^{-1} = \frac{8}{\min\left(\frac{1}{16c},2-\sqrt{2}\right)}.$$

We note that this bound is by no means an optimal one -- we saw that for a Hilbert space one can simply take $d:=1$, whereas the formula would give $d_c := 64$ (using $c:= \frac{1}{2}$, taken from the usual modulus of smoothness $\rho(\tau):=\sqrt{1+\tau^2}-1 \leq \frac{\tau^2}{2}$). Still, the above argument shows that there is a simple method one can use to readily obtain a suitable $d \geq 1$ given the original smoothness constant $c$.

\section{Acknowledgements}

The author is grateful to Ulrich Kohlenbach and Lauren\c tiu Leu\c stean for their suggestions that greatly contributed to the final form of this paper, and to the anonymous reviewer that pointed out connections to some results established in the literature.

This work was supported by a grant of the Romanian National Authority for Scientific Research, CNCS - UEFISCDI, project number PN-II-ID-PCE-2011-3-0383.

\end{document}